\documentclass[12pt]{article}
\usepackage{latexsym}
\usepackage{amssymb}
\usepackage{amsthm,amsmath}

\makeatletter
\@addtoreset{equation}{section}

\ifx\reset@font\undefined
      \let\reset@font=\relax
\fi

\title{
Restricted isometry property of matrices with independent columns
and neighborly polytopes by random sampling}

\author{
R. Adamczak \and  A. E. Litvak \and A. Pajor
       \and
N. Tomczak-Jaegermann${}^{1}$ }

\newcommand\address{\noindent\leavevmode%
\noindent
Rados{\l}aw  Adamczak, \\
Institute of Mathematics, \\
University of Warsaw, \\
Banacha 2, 02-097 Warszawa, Poland\\
 \texttt{\small
e-mail:  radamcz@mimuw.edu.pl}

\medskip
\noindent
Alexander E. Litvak, \\
Dept.~of Math.~and Stat.~Sciences,\\
University of Alberta, \\
Edmonton, Alberta, Canada, T6G 2G1.\\
\texttt{\small%
e-mail:   alexandr@math.ualberta.ca}

\medskip
\noindent
Alain  Pajor, \\
Universit\'{e} Paris-Est\\
\'{E}quipe d'Analyse et Math\'{e}matiques Appliqu\'ees, \\
5, boulevard Descartes,
Champs sur Marne,\\
77454 Marne-la-Vall\'{e}e,  Cedex 2, France\\
\texttt{\small%
e-mail: Alain.Pajor@univ-mlv.fr }

\medskip
\noindent
Nicole  Tomczak-Jaegermann, \\
Dept.~of Math.~and Stat.~Sciences,\\
University of Alberta, \\
Edmonton, Alberta, Canada, T6G 2G1.\\
\texttt{\small%
e-mail:    nicole.tomczak@ualberta.ca}

}

\date{}

\newtheorem{fact}{Fact}[section]
\newtheorem{thm}[fact]{Theorem}

\newtheorem{prop}[fact]{Proposition}
\newtheorem{lemma}[fact]{Lemma}

\newtheorem{example}[fact]{Example}
\newtheorem{definition}[fact]{Definition}

\newtheorem{Theorem}[fact]{Theorem}
\newtheorem{Lemma}[fact]{Lemma}

\newbox\nrmbox
\setbox\nrmbox=\hbox{$\Vert$}
\def\nrmrule{\vrule height\ht\nrmbox depth1.2\dp\nrmbox}

\setbox\nrmbox=%
\hbox{\kern0.15em\nrmrule\kern0.15em\nrmrule%
   \kern0.15em\nrmrule\kern0.15em}
\newcommand{\Snorm}[1]%
      {\copy\nrmbox#1\copy\nrmbox\kern-0.03em%
            \lower.4ex\hbox{}}

\renewcommand{\qed}{\bigskip\hfill\(\Box\)}

\newcommand{\R}{\mathbb R}

\newcommand{\E}{\mathbb E}
\newcommand{\EE}{B}
\newcommand{\PP}{\mathbb P}
\newcommand{\Rn}{\R^n}
\newcommand{\pp}{\mbox{{\it I\kern -0.25emP}}}

\newcommand{\eps}{\varepsilon}

\def\la{\left\langle}
\def\ra{\r\rangle}

\renewcommand{\th}{\theta}

\newcommand{\supp}{\mathop{\rm supp}}

\def\r{\right}

\def\am{A_m}
\def\bm{\EE_m}
\def\cm{C_m}

\begin{document}

\maketitle
\footnotetext[1]{This author holds the Canada Research Chair in
  Geometric Analysis.}

\begin{abstract}
This paper considers compressed sensing matrices and neighborliness
of a centrally symmetric convex polytope generated by  vectors $\pm
X_1,\dots,\pm X_N\in\R^n$, ($N\ge n$). We introduce a class of
random sampling matrices and show that they satisfy a restricted
isometry property (RIP) with overwhelming probability. In
particular, we prove that matrices with i.i.d. centered and variance
1 entries that satisfy uniformly a sub-exponential tail inequality
possess this property RIP with overwhelming probability. We show
that such ``sensing" matrices are valid for the exact reconstruction
process of $m$-sparse vectors via $\ell_1$ minimization with $m\le
Cn/\log^2 (cN/n)$. The class of sampling matrices we study includes
the case of matrices with columns that are independent isotropic
vectors with log-concave densities. We deduce that if $K\subset
\R^n$ is a convex body and $X_1,\dots, X_N\in K$ are i.i.d. random
vectors uniformly distributed on $K$, then, with overwhelming
probability, the symmetric convex hull of these points is
an $m$-centrally-neighborly
polytope with $m\sim n/\log^2 (cN/n)$.
\end{abstract}



\section{Introduction}
\label{intro}

Let $1\le m\le n\le  N$ be integers and let $X_1,\dots,X_N\in\R^n$.
Denote by  $A$ the $n\times N$ matrix with $X_1,\dots,X_N$ as
columns and by $K(A)=K(X_1,\dots,X_N)$ the convex hull of $\pm
X_1,\dots,\pm X_N$. Recall that a centrally symmetric convex
polytope is $m$-centrally-neighborly if any set of less than $m$
vertices containing no-opposite pairs, is  the vertex set of a face
(see the books \cite{Gru} and \cite{Z}).

The connection between the neighborliness of $K(A)$ and sparse
solutions of underdetermined linear equations was discovered in
\cite{D2}, Theorem 1, where it is proved that the following two
statements are equivalent:
\begin{itemize}
\item [\em i)]$K(A)$ has $2N$ vertices and is $m$-neighborly
\item [\em ii)]whenever $y=Az$ has a solution $z$
having at most $m$ non-zero coordinates
 (in other words $z$ is $m$-sparse), then
$z$ is the unique solution of the program:
$$
(P)\qquad  \min \|t\|_{\ell_1},\quad A t=A z.
$$
\end{itemize}
Here the ${\ell_1}$-norm is defined by $\|t\|_{\ell_1}=\sum_{i=1}^N
|t_i|$ for any $t=(t_i)_{i=1}^N \in\R^N$.

Statement {\em ii)} is the so-called exact reconstruction problem by
$\ell_1$ minimization or {\em basis pursuit} algorithm. For a more
detailed and complete analysis of the reconstruction of sparse
vectors by the basis pursuit algorithm we refer to
\cite{CT1} and \cite{D3}.

Let us also mention in the same stream of ideas that problem {\em
  ii)} is dual to the problem of decoding by linear programming.  In
this latter problem a linear code is given by the matrix $A^*$, and
thus a vector $x\in \R^n$ generates the vector $A^*x\in\R^N$ defined
by measurements $\bigl(\langle X_1,x\rangle,\dots, \langle
X_N,x\rangle \bigr)$. Suppose that $A^*x$ is corrupted by a noise
vector $z\in\R^N$ which is assumed to be $m$-sparse.  The problem is
to {\em reconstruct} $x$ from the data, which is the noisy output
$y=A^*x+z$.
This problem is  then  tackled by a linear
programming approach (see \cite{CT2} for complete references) that
consists of the following minimization problem
$$
(P') \qquad  \min_{t\in \R^{n}} \|y-A^*t\|_{\ell_1}.
$$

Let us denote by  $|\cdot|$ the natural Euclidean norm in $\R^n$ and
$\R^N$. Looking for a sufficient condition for a given matrix $M$ to
satisfy condition {\em ii)}, the authors of \cite{CT2} introduced
the concept of {\em Restricted Isometry Property} (RIP)   defined by
the following parameter.


\medskip
\noindent\textbf{Definition.}
{\it Let $M$ be a $n\times N$ matrix.
For any $1\le m\le \min(n,N)$, the { isometry constant} of $M$ is
defined as the smallest number $\delta_m=\delta_m(M)$ so that
$$
  (1-\delta_m)|z|^2\le {|Mz|^2}    \le( 1+\delta_m) |z|^2
$$
holds for all $m$-sparse vectors $z\in\R^N$.  The matrix $M$ is said
to satisfy the Restricted Isometry Property of order $m$ with
parameter $\delta$, shortly $\ {\rm RIP}_m (\delta)$, if
$0\leq\delta_m(M)<\delta$.
}
\medskip

The relevance of this parameter for the reconstruction property {\em
  ii)} is for instance revealed in \cite{CT1},\cite{CT2}, where it was
shown that if $\delta_{m}(M)+\delta_{2m}(M)+\delta_{3m}(M)<1$ then $M$
satisfies {\em ii)} (see also \cite{CRT}, \cite{De}, \cite{Ka}). In
the present paper, we shall use the following sufficient condition
from \cite{Candes}: if a matrix $M$ satisfies
$$
   \delta_{2m}\left({M}\right)<\sqrt 2-1
$$
then {\em i)} and {\em ii)} are satisfied.
 In other words, if $M$ has
RIP$_{2m}(\sqrt 2-1)$ then $M$ has the reconstruction property {\em
  ii)}. This approach gives the strategy of our paper.

Recall that no general construction of centrally symmetric polytopes
is known to produce polytopes with an optimal order of
neighborliness. All known results are of randomized nature,
namely, they show that for a certain probability on the space of
$n\times N$ matrices, a polytope $K(A)$ is $m$-neighborly with
overwhelming probability, for (large) $m$ depending on $n$ and
$N$. Consequently, from now on, $A$ will be a random matrix in some
{\em Ensemble} in the sense of Random Matrix Theory.  Due to the
normalization, we shall consider the isometry constant of ${A/\sqrt
  n}$. The plan consists in specializing to some model of random
matrices, the condition $\delta_{2m}\left({A\over \sqrt
    n}\right)<\sqrt 2-1$.

\medskip

Let $X_1,\dots,X_N\in \R^n$ be independent random vectors
normalized so that $\E|X_i|^2=n$ for all
$i=1,\dots, N$.
The model we will develop here is structured by two conditions: an
inequality  of the tails of linear forms and an inequality of
concentration of the Euclidean norm.
\begin{itemize}
\item  Linear forms obey a uniform sub-exponential decay,
that is, for all $1\le i\le N$, all ${y\in S^{n-1}}$, and $t>0$,
$$
   \PP\left( |\la X_i, y \ra |>t\right)\le C\exp(-ct),
$$
where $C,c>0$.

\item  The Euclidean norms of $X_1,\dots,X_N$ are
concentrated around their average:
$$
  \PP\left( \max_{i\le N}\left|{|X_i|^2\over n}-1\right|\geq{\sqrt2-1\over
   2} \right) <\lambda.
$$
\end{itemize}
Note that such a concentration inequality is clearly necessary in
order to have RIP$_{1}((\sqrt 2-1)/2)$.

One of the main results of this paper, Theorem \ref{thm:neighborly},
claims that under these conditions, whenever
$$
    m \le cn/\log^2(CN/n),
$$
the random polytope $K(A)$ is $m$-centrally-neighborly with
probability larger than $1-2\lambda -C\exp(-c\sqrt n)$, where
$C, c>0$ are universal numerical constants. We will make it more precise in
Section \ref{section:neighborly}. This model includes the cases when
\begin{itemize}
\item $X_i$'s are independent isotropic random vectors with a
  log-concave density;
\item the entries of the matrix are independent, centered with
  variance one and satisfy a sub-exponential tail inequality;
\item $X_i$'s are on the sphere of radius $\sqrt n$ and linear forms
  exhibit a uniform sub-exponential tail inequality.
\end{itemize}
These examples give rise to new classes of compressed sensing
matrices.  The class of i.i.d. entries with sub-exponential tail
behavior (that is, entries being $\psi_1$ random variables), contains
a subclass of matrices with i.i.d. $\psi_r$ entries for $1 < r \le 2$
(see Definition \ref{def:psi1_norm} below of $\psi_r$ random
variables).  Since in this case the obtained bounds are better by a
power of logarithm that may be essential in applications,
we prove our results
in full generality, for $1 \le r \le 2$.

Sub-gaussian matrices with independent $\psi_2$ entries, which correspond
to  $r=2$,
are  by now well understood. They  include for
instance
the Gaussian case when the matrix $A$ is built with i.i.d. Gaussian
$N(0,1)$ random variables (see \cite{D3},\cite{CT2},\cite{RV}); the
case when the entries of $A$ are i.i.d. $(\pm 1)$ Bernoulli random
variables (\cite{CT2}, \cite{MPT0}, \cite{BDDW}); a general case of
i.i.d. sub-gaussian entries is treated in (\cite{MPT0},\cite{MPT1},
also see \cite{MPT2} for simpler proofs).

Results of this paper are based on concentration type inequalities
for random matrices under consideration. The proof of the main
technical result, Theorem \ref{UUP_thm}, will employ methods from
\cite{ALPT}. A crucial new ingredient consists of an analysis of the
quantity
$$
  \EE_m: = \sup_{z\in U_m}\left| \left|\sum_{i \le N} z_i X_i
  \right|^2 -\sum_{i \le N} z_i^2|X_i|^2
  \right|^{1/2},
$$
 where $U_m$ denotes the set of norm one $m$-sparse vectors in $\R^N$.
 In Section~\ref{section:notation}
we present some
 definitions and preliminary tools.
 In Section \ref{section:RIP} we
 apply Theorem \ref{UUP_thm} to estimate the isometry constant
 (Theorem \ref{RIP}). Then we study the $m$-neighborly property of
 random polytopes in Section \ref{section:neighborly} and give
 application to polytopes generated by random points from a convex
 body, polytopes generated by independent vectors with independent
 $\psi_r$ random coordinates, and polytopes generated by independent
 $\psi_r$ random vectors on a sphere.
 Section \ref{section:main technical result} is devoted to the proof of
 Theorem \ref{UUP_thm} and discussion of optimality of the result.


\section{Notation and preliminaries}
\label{section:notation}

We equip $\R^n$ and $\R^N$ with the natural scalar product
$\langle\,\cdot, \,\cdot\rangle$ and the natural Euclidean norm
$|\cdot|$. We use the same notation $|\cdot|$ to denote the
cardinality of a set.
Unless otherwise stated,
$(X_i)_{i\geq 1}$ will denote independent
random vectors in $\R^n$.  By $\|M\|$ we shall denote the operator
norm of a matrix $M$, that is, $\|M\|=\sup_{|y|=1}|My|$.

\begin{definition}
\label{def:psi1_norm}
For a random variable $Y\in\R$ and $r>0$ we
define the $\psi_r$-norm by
$$
\|Y\|_{\psi_r}=\inf\left\{C>0\,;\, \,\E\exp\left({|Y|/C}\right)^r
\leq 2\right\}.
$$
\end{definition}


It is well known that the $\psi_r$-norm of a random variable may be
estimated from the growth of the moments. More precisely if a random
variable $Y$ is such that for any $p\geq 1$, $\|Y\|_p\le p^{1/r} K$,
for some $K >0$, then $\|Y\|_{\psi_r}\le c K$ where $c>0$ is a
numerical constant.

\begin{definition}
\label{def:psi_vector norm}
Let $X\in\R^n$ be a centered random vector and $r >0$. We say that
  $X$ is  $\psi_r$ or a $\psi_r$ vector, if
 $ \ \sup_{y\in S^{n-1}} \|\la X, y \ra \|_{\psi_r}$
is bounded and we set
$$
\|X\|_{\psi_r}=\sup_{y\in S^{n-1}} \|\la X, y \ra \|_{\psi_r}.
$$
\end{definition}


\noindent {\bf Remark:}
The above notation of $\|X\|_{\psi_r}$ for the {\em weak} $\,\psi_r$
{\em norm} of a random vector $X$ should not be confused with the
standard convention in the probability theory that this notation
stands for the $\psi_r$ norm of the random variable $|X|$, i.e.,
$\|\,|X|\,\|_{\psi_r}$--this latter meaning will never be used in this
paper.

\bigskip
We recall the well known Bernstein's inequality which we shall use in
the form of a $\psi_1$ estimate (\cite{vVW}).

\begin{Lemma}
    \label{lemma:Bernstein}
Let $Y_1,...,Y_n$ be independent real random variables with zero mean
such that for some $\psi >0$ and every $i$, $\|Y_i\|_{\psi_1} \le \psi$.
Then, for any  $t>0$,
\begin{equation*}
\PP \Bigl(\bigl|\sum_{i=i}^n Y_i \bigr| >t
\Bigr) \leq 2\exp \left(-\,
\frac{t^2}{4\sum_{i\le n}\|Y_i\|_{\psi_1}^2+2t\psi}\right).
\end{equation*}
\end{Lemma}

Given a set $E\subset \{1, ..., N\}$ by $P_E$ we denote the orthogonal
projection from $\R^N$ onto  the coordinate subspace of vectors whose
supports  are  in $E$. We denote this  subspace by $\R^E$.
The support of  $z\in\R^N$  is  denoted by $ \supp  z$.
A vector
$z \in \R^N$  is  called $m$-sparse
if $|\supp z|\le m$. The subset of $m$-sparse unit vectors in  $\R^N$
is denoted by
\begin{equation}\label{U_m}
U_m=U_{m}(\R^N): =\{z\in\R^N\,:\, |z|=1, | \supp z|\le m\}.
\end{equation}

Let $B_\infty^N=\{x=(x_i)\in\R^N \,:\,\|x\|_\infty=\max_i|x_i|\le 1\}$
and $B_2^N$ be the unit Euclidean ball.  For every $E\subset \{1, ...,
N\}$, $\eps, \alpha \in (0, 1]$ by ${\cal{N}}(E, \eps, \alpha)$ we
denote an $\eps$-net in $B_2^N \cap \alpha B_{\infty}^N \cap \R^E$ in
the Euclidean metric. Thus for every $x\in B_2^N \cap \alpha
B_{\infty}^N$ supported by $E$, there exist $\bar x\in {\cal{N}}(E,
\eps, \alpha)$ supported by $E$ such that $|x-\bar x|<\eps$.  A
standard volume comparison argument shows that we may assume that the
cardinality of ${\cal{N}}(E, \eps, \alpha)$ does not exceed
$(3/\eps)^m$, where $m$ is the cardinality of $E$.

\begin{definition}
\label{isotropic_vector}
A random vector $X\in \R^{n}$ is called isotropic if
\begin{equation}
\E\langle X,y\rangle=0,\quad \E\,|\langle X, y \rangle|^{2}=|y|^{2}
\quad \mbox{\rm for all }
 y\in \R^{n},  \label{isotropic}
\end{equation}%
in other words, if $X$ is centered and its covariance matrix is the
identity.
\end {definition}
A subset $K\subset\R^n$ is said to be isotropic when  a  random
point $X$ uniformly distributed in $K$ is an isotropic random
vector.

Recall that a function $f: \R^n \to \R $ is
called log-concave if for any $\theta \in \lbrack 0,1]$ and any
$x_{1},x_{2}\in \mathbb{R}^{n}$,
\begin{equation*}
f\big(\theta x_{1}+(1-\theta )x_{2}\big)\geq f(x_{1})^{\theta
}f(x_{2})^{1-\theta }.
\end{equation*}

It is well known \cite{Bor2} that if a measure has a log-concave
density, then linear functionals exhibit a sub-exponential decay.
More precisely, we have:

\begin{lemma}\cite{Bor2}:\label{psi-lemma}
  Let $X\in\R^n$ be a centered random vector with a log-concave
  density.  Then for every $y \in S^{n-1}$,
$$
  \|\la X, y \ra \|_{\psi_1} \leq c\,\left(\E|\langle X, y \rangle
  |^2\right)^{1/2},
$$
where $c>0$ is a universal constant. As a consequence, if $X$ is an
isotropic random vector with a log-concave density then
$\|X\|_{\psi_1}\le c$.
\end{lemma}

The Euclidean norm of an isotropic random vector with a log-concave
density highly concentrates around its expectation, this translates
geometrically to the concentration of mass of an isotropic convex
body within a thin Euclidean shell (\cite{Kla2}, see also
\cite{FGP}). We will use here the following result  immediately
derived from \cite{Kla1}, Theorem 4.4.

\begin{lemma}\label{kla}
  Let $1\le n\le N$ be integers and let $X_1,\dots, X_N\in\R^n$ be
  isotropic random vectors with log-concave densities. There
  exist numerical positive constants $C,c_0$ and $c_1\in (0,{1\over
    2})$ such that for all $\theta\in(0,1)$ and $ N\le
  \exp(c\theta^{c_0} n^{c_1})$,
$$
 \PP\left( \max_{i\le N} \left| {|X_i|^2\over n} -1\right|\ge
 \theta\right) \le C\exp (-c\theta^{c_0}n^{c_1}).
$$
Moreover, one can take $c_0=3.33$ and $c_1=0.33$.
\end{lemma}

\noindent {\bf Remark:} It is conjectured that in the above theorem
one can replace $\theta^{3.33}n^{0.33}$ by $c(\theta)n^{1/2}$.

\medskip

We shall also use the following result from \cite{Pao} as formulated
in \cite{ALPT}.

\begin{lemma}
\label{max}
 Let $N,n\ge 1$ be integers and let $X_1,\dots, X_N \in \R^n$
 be isotropic random vectors with log-concave densities.
 Then there exists an absolute positive constant $C_0$ such that for any
  $N\leq \exp(\sqrt{n})$ and for every $K\geq 1$ one has
$$
   \max _{i\leq N} |X_i| \leq C_0 K \sqrt{n}
$$
with probability at least $1 - \exp(-K \sqrt{n})$.
\end{lemma}

In this paper, different universal positive constants may be denoted by the
same letters $C,C_0,C', c,c_0,c'$, etc.

\section{Isometry constant}
\label{section:RIP}

We begin this section by formulating, in Theorem \ref{RIP}, a general
estimate for the isometry constant of random matrices with independent
$\psi_r$ columns. Then, in order to apply such an estimate, we
introduce two sufficient conditions that determine large classes of
random matrices.  Finally, we give examples of important classes that
satisfy the estimates from Theorem \ref{RIP} and thus provide us with
models: the Log-Concave Ensemble, matrices with i.i.d. $\psi_r$
entries, and matrices defined by independent $\psi_r$ vectors on a
sphere.

\subsection{Estimating the isometry constant}

Techniques of ``compressed sensing" rely on properties of the
sampling matrix, which should act almost isometrically on sparse
vectors.  This motivated the concept of Restricted Isometry Property
(RIP) defined in \cite{CT2}. To quantify this property of the
``sensing" matrix, the authors introduced the isometry constant
defined in the introduction, that we recall here  for the convenience
of the reader.

\begin{definition}
\label{rip-def} Let $M$ be a $n\times N$ matrices and let $\delta\in
(0,1)$. For any $1\le m\le \min(n,N)$, the { isometry constant} of
$M$ is defined as the smallest number $\delta_m=\delta_m(M)$ so that
\begin{equation}
\label{equ:RIP} (1-\delta_m)|z|^2\le {|Mz|^2}    \le( 1+\delta_m)
|z|^2
\end{equation}
holds for all $m$-sparse vectors $z\in\R^N$.  The matrix $M$ is said
to satisfy the Restricted Isometry Property of order $m$ with
parameter $\delta$ if $0\leq\delta_m(M)<\delta$.
\end{definition}

Let $X_1,\ldots,X_N\in \R^n$ and let $A=A^{(n,N)}$ be the
``sampling" matrix with the $X_i$'s as columns.  We begin by a
simple observation. Define the following quantity
\begin{align}\label{B_m_definition}
  \EE_m= \sup_{z\in U_m}\left| \left|\sum_{i\le N}z_iX_i
  \right|^2-\sum_{i\le  N}z_i^2|X_i|^2
 \right| ^{1/2}.
\end{align}
 Then, clearly
\begin{equation}
\label{equ: RIP from Bm}
\sup_{z\in U_m}\left| {|Az|^2\over n}-1
 \right|\le {\EE_m^2\over n}+ \max_{i \le  N} \left|{|X_i|^2\over n}
   -1\right|.
\end{equation}
Thus the isometry  constant  is controlled by quantity  $B_m$
and the second term,  $\max_{i \le  N} \left|{|X_i|^2\over n}
   -1\right|.$ We begin by estimating   $\EE_m$.

\begin{Theorem}
\label{UUP_thm}
Let $n\geq 1$ and $1\le m\leq N$ be integers. Let
$1\leq r \leq 2$ and $X_1,\ldots,X_N\in \R^n$ be independent
$\psi_r$ random vectors with $\psi=\max_{i \le N}\|X_i\|_{\psi_r}$.
Let $\theta\in (0,1/4)$,  $K,K'\geq 1$  and  assume that  $m$ satisfies
$$
   m\log^{2/r} \frac{2N}{\theta m}\le \theta^2 n.
$$
Then setting   $\xi=\psi K +K'$, the inequality
$$
 \EE_m ^2\leq  C\xi^2 \theta n
$$
holds with probability at least
$$
  1-  \exp\left( - c K^r \sqrt{m}\log^{} \left(\frac{2N}{\theta
        m}\right) \r)-
  \PP \left( \max_{i\le N} |X_i|\ge K'\sqrt n\right),
$$
where $C,c$ are absolute positive constants.
\end{Theorem}

We postpone the proof of Theorem \ref{UUP_thm} to the last section.
Combining this theorem with inequality (\ref{equ: RIP from Bm}),
relating the RIP, $B_m$ and concentration of the Euclidean norm of the
$X_i$'s, we  immediately  deduce   an estimate for the
isometry constant of a random matrix with independent $\psi_r$
columns.

\begin{thm}
 \label{RIP}
 Let $n\geq 1$ and $m,N$ be integers such that $1 \le m \le
 \min(N,n)$. Let $1\le r\le 2$. Let $X_1,\ldots,X_N\in \R^n$ be
 independent $\psi_r$ random vectors and let $\psi=\max_{i\le
   N}\|X_i\|_{\psi_r}$.  Let $\theta'\in (0,1)$, $K,K'\geq 1$ and set
 $\xi=\psi K +K'$.  Then
 $$
 \delta_m\left({A\over\sqrt n}\right)\le C\xi^2 \sqrt{m\over
   n}\log^{1/r}\left({eN\over  m\sqrt{m\over n}}\right)  + \theta'
$$
holds with probability larger than
\begin{eqnarray*}
  1 &-& C\exp\left( - cK ^r\sqrt m \log^{}\left(\frac{eN}
          {m\sqrt{\frac{m}{ n}}}\right) \r) \\
&-&    \PP \left( \max_{i\le  N} |X_i|\ge K'\sqrt n\right)
  -\PP \left( \max_{i\le  N} \left|{|X_i|^2\over n} -1\right|\ge
  \theta'\right),
\end{eqnarray*}
where $C,c>0$ are universal constants.
\end {thm}

Note that such an estimate for $ \delta_m\left({A\over\sqrt n}\right)$
is meaningful only if, firstly, its right hand side is $<1$, and
secondly, if it holds with probability $>0$.  In fact, the former
condition is equivalent to the RIP of order $m$.  This leads to
considerations of models of random $n\times N$ matrices that satisfy
the following two conditions.  Let $1\le r\le 2$, $\psi>0$ and
$\lambda\in (0,1)$. Let $X_1,\ldots,X_N\in\R^n$ be independent
$\psi_r$ random vectors and let $A$ be the matrix with
$X_1,\ldots,X_N$ as columns.
\begin{itemize}
\item {\bf Condition $H_1(r,\psi)$:}
Linear forms obey a uniform $\psi_r$ estimate:
\begin{equation} \label{equa:sub-exp hyp}
\|\langle X_i,y\rangle \|_{\psi_r}\le \psi \quad \text{for all}
\quad y\in S^{n-1} \ \text{and}\ 1\le i\le N.
\end{equation}
\item {\bf Condition $H_2(\lambda)$:} $|X_i|$'s are concentrated
  around their average:
\begin{equation} \label{equa:concentration hyp}
\PP\left( \max_{i\le N}\left|{|X_i|^2\over
n}-1\right|\geq{\sqrt2-1\over 2} \right) <\lambda.
\end{equation}
\end{itemize}

As already mentioned in the Introduction, a condition such as
$H_2(\lambda)$ is  necessary to have the RIP.  Indeed, if
the matrix $A/\sqrt n$ has RIP$_{1}((\sqrt 2-1)/2)$ with probability
$\lambda$ then $H_2(\lambda)$ is satisfied.

\subsection{Examples}

We now specialize Theorem  \ref{RIP} to some specific classes of
matrices.

\subsubsection{The Log-Concave Ensemble}

We start by considering the ``log-concave  setting'', where $X_1,\dots,
X_N \in \R^n$ are independent isotropic vectors with log-concave
densities.

\begin{lemma}
  \label{logconcave:setting}
  Assume the above ``log-concave setting''.  There exist universal
  constants $\psi, C, c>0$ such that conditions $H_1(1,\psi)$ and
  $H_2(C\exp(-cn^{c_1}))$ are satisfied  whenever  $N\le
  \exp{(cn^{c_1})}$, where ${c_1}$ is given in Lemma \ref{kla}.
\end{lemma}

The  proof  is  immediate from Lemmas \ref{psi-lemma}  and \ref{kla}.

\medskip

Applying Theorem \ref{RIP} (with $r=1$)  together with Lemmas
\ref{logconcave:setting} and \ref{max}
to the Log-Concave Ensemble, we get that
for every $N\le \exp(cn^{c_1})$,
\begin{equation}\label{logconcave:RIP}
\delta_m\left({A\over \sqrt n}\right)\le C \sqrt{m\over
  n}\log\left({eN\over
    m\sqrt{m\over  n}}\right) + {\sqrt2-1\over 2}
\end{equation}
holds with probability larger than
$$
  1- C\exp\left( - c \sqrt{m}\log\left({eN\over m\sqrt{m\over
  n}}\right) \r)-  e^{-c\sqrt n} -  \exp(-cn^{c_1}),
$$
where $C, c>0$ are universal constants and $c_1$ is given in
Lemma \ref{kla}.

It might be worthwhile to note that using directly Lemma \ref{kla} one
can replace the second term in estimate (\ref{logconcave:RIP}) by a
term tending to 0 when $n \to \infty$, but this would require an
adjustment in probability. For example $1/n^{c_1/2c_0}$ works with the
probability estimate in which $\exp(-cn^{c_1})$ is replaced by
$\exp(-cn^{c_1/2})$. (Here $c_0$ is given in Lemma \ref{kla}.)

\subsubsection{Matrices  with independent $\psi_r$ entries}
\label{iid-entries-rip}

Consider now the ``$\psi_r$  setting", where the entries $a_{ij}$ of the
matrix $A$ are independent centered, with variance one, random
$\psi_r$ variables (with $ r \in [1,2]$).  Set $\psi=\max_{ij}
\|a_{ij}\|_{\psi_r}$.

\begin{lemma}
\label{psi1:setting}
Assume the above ``$\psi_r$  setting'' with $ r \in [1,2]$.
Then
conditions $H_1(r,C\psi)$ and $H_2(2\exp(-c n^{r/2}/\psi^{2 r}))$
are satisfied   whenever  $N\le \exp(cn^{r/2}/\psi^{2r})$, where
$C,c$ are absolute positive constants.
\end{lemma}

\noindent{\bf Proof.}
To prove that the columns of the matrix $A$ are $\psi_r$ vectors we
will estimate the $p$-th moments of random variables $\sum_{i=1}^n y_i
a_{ij}$, for any $y=(y_i) \in \R^n$ and any $p \ge 1$. This will be
done by using Talagrand's concentration inequality for linear
combinations of symmetric Weibull variables together with some
symmetrization and truncation arguments.

The following Lemma is a combination of Corollaries 2.9 and 2.10 of \cite{T2}.

\begin{lemma}\label{Weibull}
Let $r \in [1,2]$ and $Y_1,\ldots,Y_n$ be independent symmetric
random variables satisfying $\PP(|Y_i| \ge t) = \exp(-t^r)$. Then
for every vector $a = (a_1,\ldots,a_n) \in \R^n$ and every $t \ge 0$,
\begin{displaymath}
 \PP\Big(\Big|\sum_{i=1}^n a_i Y_i\Big| \ge t \Big) \le
 2\exp\Big(-c\min\Big(\frac{t^2}{\|a\|_2^2},
 \frac{t^r}{\|a\|_{r^\ast}^r}\Big)\Big),
\end{displaymath}
where $1/r^\ast + 1/r = 1$ and $\|a\|_q = (|a_1|^q + \ldots
+|a_n|^q)^{1/q}$, for $1\le q <\infty$.
\end{lemma}

The behavior of general centered $\psi_r$ variables can be easily
reduced to symmetric Weibull variables. The argument is quite
standard, we sketch it here for the sake of completeness.

Assume thus that $Z_1,\ldots,Z_n$ are independent mean zero random
variables with $\|Z_i\|_{\psi_r}\le 1$. Let $\beta = (\log 2)^{1/r}$
and   set  $U_i = (|Z_i| - \beta)_+$. Let $Y_i$ be defined as in Lemma
\ref{Weibull}.

We have for $t > 0$,
\begin{align*}
\PP(U_i \ge t) &\le \PP(|Z_i| \ge t + \beta) \le 2\exp(-(t+\beta)^r) \\
&\le 2\exp(-t^r - \beta^r) = \PP(|Y_i|\ge t).
\end{align*}

We will use the above observation together with symmetrization and the
contraction principle to estimate moments of linear combinations of
variables $Z_i$. We have for $p \ge 1$,
\begin{align*}
  \Big\|\sum_{i=1}^n a_i Z_i \Big\|_p &\le 2\Big\|\sum_{i=1}^n a_i
  \varepsilon_i Z_i \Big\|_p \quad \textrm{(symmetrization)}\\
  & = 2\Big\|\sum_{i=1}^n a_i \varepsilon_i |Z_i| \Big\|_p \\
  & \le 2\Big\|\sum_{i=1}^n a_i \varepsilon_i (\beta + U_i) \Big\|_p
  \quad \textrm{(the contraction principle})\\
  & \le 2\Big\|\sum_{i=1}^n a_i \varepsilon_i \beta \Big\|_p +
  2\Big\|\sum_{i=1}^n a_i \varepsilon_i U_i \Big\|_p \\
  & \le C\sqrt{p}\beta \|a\|_2 + 2\Big\|\sum_{i=1}^n a_i \varepsilon_i
  Y_i \Big\|_p \\
  & \le C\sqrt{p}\|a\|_2 + Cp^{1/r}\|a\|_{r^\ast},
\end{align*}
where to get the last two inequalities we used Khinchine's inequality,
Lemma \ref{Weibull} and integration by parts to pass from tail to
moment estimates.

We are now ready to prove condition $H_1(r,C\psi)$. Fix $y \in
S^{n-1}$ and consider the linear combination $\sum_{i=1}^n y_i
a_{ij}$.  Since $\|a_{ij}\|_{\psi_r} \le \psi$,  we
obtain  by homogeneity
\begin{displaymath}
\Big\|\sum_{i=1}^n y_i a_{ij}\Big\|_p \le C\psi(\sqrt{p}\|y\|_2 +
Cp^{1/r}\|y\|_{r^{\ast}}) \le 2C\psi p^{1/r},
\end{displaymath}
since $r \in [1,2]$ implies that $p^{1/r} \ge \sqrt{p}$ and
$\|a\|_{r^{\ast}} \le \|a\|_2 = 1$. The growth condition on the
moments of the random variable $\sum_{i=1}^n y_ia_{ij}$ implies that
its $\psi_r$ norm is bounded by $\tilde{C}\psi$.  \smallskip

The proof of condition $H_2$ goes along similar lines. Instead of
Lemma~\ref{Weibull} we will now use the following lemma, which is an
easy consequence of Theorem 6.2 in \cite{H} and the observation that
the $p$-th moment of a Weibull variable with parameter $s$ is of order
$C_s p^{1/s}$, where $C_s$ remains bounded for $s$ away from $0$.

\begin{lemma}\label{HMSO}
  If $0 < s <1$ and $Y_1,\ldots,Y_n$ are independent symmetric random
  variables satisfying $\PP(|Y_i| \ge t) = \exp(-t^s)$, then for $a =
  (a_1,\ldots,a_n) \in \R^n$ and $p \ge 2$,
\begin{displaymath}
\Big\|\sum_{i=1}^n a_i Y_i \Big\|_p \le C\sqrt{p}\|a\|_2 + C_s
p^{1/s}\|a\|_{p}.
\end{displaymath}
Moreover, for $s \ge 1/2$, $C_s$ is bounded by some absolute constant.
\end{lemma}

Using similar arguments as in the proof of condition $H_1$ we can
infer from the above lemma that if $Z_1,\ldots,Z_n$ are independent
mean zero random variables with $\|Z_i\|_{\psi_s} \le b$ ($s \in
[1/2,1)$), then for $p \ge 2$,
\begin{displaymath}
\Big\|\sum_{i=1}^n a_i Z_i \Big\|_p \le Cb(\sqrt{p}\|a\|_2 +
p^{1/s}\|a\|_{p}).
\end{displaymath}
Therefore, for any $p \ge 2$ by the Chebyshev inequality in $L_p$,
\begin{displaymath}
\PP\Big(\Big|\sum_{i=1}^n Z_i\Big| \ge Cb(\sqrt{np} + p^{1/s}
n^{1/p}) \Big) \le \exp(-p).
\end{displaymath}
For $p \ge 3$ we have
\begin{displaymath}
\sqrt{np} + p^{1/s} n^{1/p} \le \tilde{C}(\sqrt{np} + p^{1/s})
\end{displaymath}
with $\tilde{C}$ universal for $s \ge 1/2$, so the above inequality yields
\begin{displaymath}
\PP\Big(\Big|\sum_{i=1}^n Z_i\Big| \ge Cb(\sqrt{np} + p^{1/s}) \Big) \le e^3\exp(-p)
\end{displaymath}
for some (new) universal constant $C$ or equivalently
\begin{align}\label{Berstein_psi_r}
\PP\Big(\Big|\sum_{i=1}^n Z_i\Big| \ge t \Big) \le
2\exp\Big(-c\min\Big[\frac{t^2}{nb^2},
       \Big(\frac{t}{b}\Big)^{s}\Big]\Big).
\end{align}

For fixed $j$ we apply this inequality with $s = r/2$ to variables
$Z_i = a_{ij}^2 - 1$. Note that $\E Z_i = 0$ and
\begin{align*}
  \|Z_i\|_{\psi_{r/2}} &\le C(1 +
\|a_{ij}^2\|_{\psi_{r/2}}) \\
  &= C(1+ \|a_{ij}\|_{\psi_r}^2) \le \tilde{C}\psi^2.
\end{align*}
(The additional constants appearing above stem from the fact that
under the standard definition for $s < 1$, $\|\cdot\|_{\psi_{s}}$ is
not a norm but only a quasi-norm and additionally $\|1\|_{\psi_{r/2}}
\neq 1$.  One can modify the function $x \mapsto e^{x^r} - 1$ so that
it is convex.  For $r$ away from zero, this modification changes the
norm by an absolute constant).  Therefore, applying
(\ref{Berstein_psi_r}) with $t = \varepsilon n$ yields
\begin{align*}
  \PP\Big(\Big|\frac{1}{n}\sum_{i=1}^n a_{ij}^2 - 1\Big| \ge
  \varepsilon \Big) &\le 2\exp\Big(-c\min\Big[\frac{\varepsilon^2
    n}{\psi^4},\Big(\frac{\varepsilon n}{\psi^2}\Big)^{r/2}\Big]\Big)
  \\ &\le 2\exp\Big(-\tilde{c}\frac{\varepsilon^r
    n^{r/2}}{\psi^{2r}}\Big).
\end{align*}

For $r = 2$ the proof is similar, but uses Lemma \ref{Weibull} (which
in this case reduces to Bernstein's $\psi_1$ inequality) instead of
Lemma \ref{HMSO} (the argument is simpler since in this case the
involved norms of the vector $a$ do not depend on $p$ and we get
(\ref{Berstein_psi_r}) directly).

The lemma follows now by the union bound.
\qed

Applying Theorem \ref{RIP} together with Lemma \ref{psi1:setting}
to the ``$\psi_r$ setting'',
we get that for every $N\le \exp(cn^{r/2}/\psi^{2r})$,
\begin{equation}
\label{psi1:RIP}
\delta_m\left({A\over \sqrt n}\right)\le C \psi^2 \sqrt{m\over
  n}\log^{1/r}\left({eN\over
    m\sqrt{m\over  n}}\right)+{\sqrt2 -1\over 2}
\end{equation}
holds  with probability at least
$$1- C\exp(-cn^{r/2}/\psi^{2r})$$
where $C,c>0$ are universal constants.

\subsubsection{Vectors on a sphere}

Another interesting case is when the vectors $X_1,\dots, X_N$ lie on a
common sphere.  To keep the same normalization as in the previous
cases we assume that the sphere has the radius $\sqrt n$.  Then
condition (\ref{equa:concentration hyp}) becomes empty.  Let $1\le
r\le 2$ and assume that the vectors are $\psi_r$ and let $\psi=\max_{i
  \le N}\|X_i\|_{\psi_r}$. Let $K\geq 1$ and set $\xi=\psi K$.  Then
Theorem \ref{RIP} immediately gives that
\begin{equation}
\label{samenorm:RIP}
\delta_m\left({A\over \sqrt n}\right)\le C\xi^2 \sqrt{m\over
  n}\log^{1/r}\left({eN\over     m\sqrt{m\over  n}}\right)
\end{equation}
with probability larger than
$$
1- C\exp\left( - cK^r \sqrt{m}\log^{}\left({eN\over m\sqrt{m\over
        n}}\right) \r)
$$
where $C,c>0$ are universal constants.

\section{The geometry of faces of random polytopes}
\label{section:neighborly}

In this Section we discuss the geometry of random polytopes. Let
$A $ be an $n\times N$ matrix. We denote by $K^+(A)$ (resp. $K(A)$)
the convex hull (resp., the symmetric convex hull) of the $N$
columns of $A$.

\subsection{Neighborly polytopes}

For an integer  $1 \le m \le n $, a polytope is called {\em
$m$-neighborly} if any set of less than $m$ vertices is the vertex
set of a face.  In the symmetric setting, a centrally symmetric
convex polytope is $m$-centrally-neighborly if any set of less
than $m$ vertices containing no-opposite pairs is  the vertex set
of a face.
We refer the reader to the books \cite{Gru} and \cite{Z} for
classical details on neighborly  polytopes. (Some  new quantitative
invariants   related to neighborliness were recently developed in
\cite{Man_tom}.)

The relation between the problem of reconstruction and neighborly
polytopes was discovered in \cite{D2}.

\begin{Theorem} (\cite{D2}, Theorem 1)
  \label{thm:Donoho}
  Let $A$ be a $n\times N$ matrix, $n\le N$. The following two assertions
  are equivalent.
\begin{itemize}
\item [\em i)]The polytope $K(A)$ has $2N$ vertices and is
  $m$-centrally-neighborly.
\item [\em ii)]Whenever $y=Az$ has a solution $z$ having at most $m$
  non-zero coordinates, $z$ is the unique solution of the optimization
  problem $(P)$:
$$
(P)\qquad  \min \|t\|_{\ell_1},\quad A t=A z,
$$
\end{itemize}
\end{Theorem}

We will also use the following result from \cite{Candes} (which  could
be replaced by a similar result from \cite{CT2}).

\begin{lemma} \cite{Candes}
\label{candes:lemma}
 Assume that $\delta_{2m}(A/\sqrt n)<\sqrt 2-1$. Then whenever $y=Az$
 has a solution $z$ having at most $m$ non-zero coordinates, $z$ is
 the unique solution of the $\ell_1$ minimization problem $(P)$.
\end{lemma}

We are now ready to state the main result on neighborly random
polytopes.

\begin{Theorem}
  \label{thm:neighborly}
  Let $1\leq m \le n\leq N$ be integers. Let $1\le r\le 2$. Let
  $\psi \ge 1$ and $\lambda\in (0,1/2)$.
  Let $X_1,\dots, X_N$ be independent  random
  vectors satisfying $H_1(r,\psi)$ with parameter $\psi$ and
  $H_2(\lambda)$ with probability $\lambda$.  Let $A$ be the $n\times
  N$ matrix with $X_1,\dots, X_N$ as columns.  Then, with probability
  larger than
  $$
     1-2\lambda-C\exp(-c\sqrt n/\psi^2)
  $$
  the polytopes $K^+(A)$ and $K(A)$ are $m$-neighborly and
  $m$-centrally-neighborly, respectively, whenever
$$
  m\le {cn\big/ \psi^{4}\log^{2/r} (C\psi^6 N /n)},
$$
where $C, c>0$ are universal constants.
\end{Theorem}

\noindent Observe that the probability is positive for $n$ large
enough provided  that  $\lambda <1/2$.

\medskip

\noindent{\bf Proof.}
 Theorem \ref{RIP}
and the definition of property $H_1(r, \psi)$  imply that
for arbitrary $\theta' \in (0,1)$,  and  $K, K' \ge 1$,
 setting $\xi= \psi K + K'$,
we have
$$
 \delta_m\left({A\over\sqrt n}\right)\le C\xi^2 \sqrt{m\over
   n}\log^{1/r}\left({eN\over  m\sqrt{m\over n}}\right)  + \theta'
$$
holds with probability larger than
\begin{eqnarray}
  1 &-& C\exp\left( - cK^r \sqrt m \log^{}\left(\frac{eN}
          {m\sqrt{\frac{m}{ n}}}\right) \r) \nonumber\\
\label{prob_est}
&-&    \PP \left( \max_{i\le  N} |X_i|\ge K'\sqrt n\right)
  -\PP \left( \max_{i\le  N} \left|{|X_i|^2\over n} -1\right|\ge
  \theta'\right).
\end{eqnarray}
In view of Lemma \ref{candes:lemma}, we look for $m$ and $\theta'$ to
ensure $\delta_{2m}(A/\sqrt n)<\sqrt 2-1$. For instance, we let
$\theta'=(\sqrt 2-1)/2$ and note that (\ref{equa:concentration hyp})
implies
\begin{equation}
\label{proba:max}
\PP\left( \max_{i\le  N}|X_i|\geq\left({\sqrt2-1\over
      2}+1\right)^{1/2} \sqrt n \right) <{\lambda}.
\end{equation}
So we take $K'=\left({\sqrt2-1\over 2}+1\right)^{1/2}$ and $K=1$ which
determines $\xi=\psi K+K'$ in terms of $\psi$.  We shall use the fact
that $ 1 \le \xi/\psi \le \tilde C$,
where  $\tilde C$ is a universal constant.

\medskip

Now set $m_0 = [{c'n\big/ \psi^{4}\log^{2/r} (C'\psi^6 N /n)}]$ (for
some new constants $C',c'>0$). It is clearly sufficient to prove that
the polytopes $K^+(A)$ and $K(A)$ are $m_0$-neighborly and
$m_0$-centrally-neighborly, respectively.  Thus adjusting the
constants  $C',c'>0$ and writing  $m $ for $ m_0$, we obtain
$$
C\xi^2 \sqrt{m\over n}\log^{1/r}\left({eN\over  m\sqrt{m\over n}}\right)
<  (\sqrt 2-1)/2.
$$
Combining this with the choice of $\theta'$,
passing from $m$ to $2m$ and adjusting the constants again if
necessary,
 we conclude  that
$\delta_m\left({A\over\sqrt n}\right) <  \sqrt 2-1$ with probability
larger than
$$
  1- C''\exp\left( - c''\sqrt m \log^{}\left(\frac{eN}
  {m\sqrt{\frac{m}{ n}}}\right) \right) -2\lambda.
$$
The last estimate follows from (\ref{prob_est}) by applying
(\ref{equa:concentration hyp}) and (\ref{proba:max}) to the last two
terms, respectively; and where $C'',c''>0$ are again new
constants.~\qed

\subsection{Examples}

We will now apply Theorem \ref{thm:neighborly} in  the three
different settings introduced in the previous section.

\subsubsection{The Log-Concave Ensemble}

Applying Lemma \ref{logconcave:setting} and  bound
(\ref{logconcave:RIP})  we get the following:
\begin{Theorem}
  \label{thm:neighborly2}
  Let $1\le m \le n\leq N$ be integers.
  Let $X_1,\dots, X_N$ be
  independent isotropic vectors with log-concave densities.  This is
  for instance the case if $X_1,\dots, X_N$ are i.i.d. random vectors
  uniformly distributed on an isotropic convex body.  Then, for any
  $N\le \exp(cn^{c_1/2})$, with probability at least $1-
  C\exp(-cn^{c_1/2})$, the polytopes $K^+(A)$ and $K(A)$ are
  $m$-neighborly and $m$-centrally-neighborly, respectively, whenever
$$
   m\le {cn\big/ \log^2 (CN /n)},
$$
where $C,c>0$ are universal constants and ${c_1}$ is given in Lemma
\ref{kla}.
\end{Theorem}

\medskip

\noindent{\bf Remark:}  It is known (\cite{BGMN}) that there
is a universal constant $\psi$ such that the uniform probability
measure on the ball $\{ x\in\R^n\,:\, \sum_{i=1}^n |x_i|^r\leq 1\}$
satisfies $H_1(r,\psi)$ for $1\le r\leq 2$ and satisfies
$H(2,\psi)$ for $r\ge 2$. Of course, since it is log-concave, the
concentration property $H_2$ is also satisfied. Applying Theorem
\ref{thm:neighborly}
to these examples, we get a better estimate of  the level of
neighborliness than  in Theorem \ref{thm:neighborly2}. We get now
$m\sim {cn\big/ \log^{2/r} (CN /n)}$ for $1\leq r\leq 2$ and $m\sim
{cn\big/ \log (CN /n)}$
for $2\le r\le \infty$.

\subsubsection{Matrices  with independent $\psi_r$ entries}

In a similar way as above, Lemma  \ref{psi1:setting} and bound
(\ref{psi1:RIP})  imply the following theorem (note that its
conclusion becomes empty if $N\ge \exp(cn^{r/2}/\psi^{2r})$ and
$\psi \geq 1$).

\begin{Theorem}
  \label{thm:neighborly4}
  Let $A$ be a matrix with entries  that are independent
  centered  variance one random variables. Let $1\le r\le 2$ and
  assume that the $\psi_r$
  norms of the entries are  bounded by some constant $\psi$.
  Then, for any $N\le \exp(cn^{r/2}/\psi^{2r})$, with probability
  at least $1- C\exp(-cn^{r/2}/\psi^{2r})$, the polytopes
  $K^+(A)$ and $K(A)$ are $m$-neighborly and $m$-centrally-neighborly,
  respectively,   whenever $1 \le m\leq n$ satisfies
$$
  m\le {cn\big/ \psi^4\log^{2/r} (C\psi^6 N /n)} ,
$$
where $C,c>0$ are universal constants.
\end{Theorem}

\subsubsection{Vectors on a sphere}

Finally assume that the vectors are on a sphere of radius $\sqrt n$.
  From bound (\ref{samenorm:RIP}) we obtain:

\begin{Theorem}
  \label{thm:neighborly3}
  Let $1\le m \le n\leq N$ be integers. Let $1\le r\le 2$. Let $X_1,\dots,
  X_N$ be independent vectors on a sphere of radius $\sqrt n$ and
  satisfying $H_1(r,\psi)$ for some parameter $\psi>0$. Let $K\ge 1$
  and set $\xi=\psi K$.  Then, with probability at least $1-
  C\exp(-K\sqrt n/\psi^2)$, the polytopes $K^+(A)$ and $K(A)$ are
  $m$-neighborly and $m$-centrally-neighborly, respectively, whenever
$$
  m\le {cn\big/ \xi^4\log^{2/r} (C\xi^6 N /n)},
$$
where $C, c>0$ are universal constants.
\end{Theorem}

\noindent{\bf Remark:} {\bf 1)} For the matrix $A$ with i.i.d.
Gaussian $N(0,1)$ entries (the case considered in Section
\ref{iid-entries-rip} above when $r=2$), it is known that with
overwhelming probability, $K(A)$ is $m$-centrally-neighborly, whenever
$1 \le m\leq n$ satisfies
$$
  m\le {cn\big/ \log (CN /n)},
$$
where $C,c>0$ are universal constants, (see
\cite{D3},\cite{CT2},\cite{RV},\cite{MPT1}). The precise asymptotic
dependence of $m$ on $n$ and $N$ has been well studied in \cite{D4} when
$n/N\to\delta\in(0,1)$ and in \cite{DT} when $n/N\to 0$.

{\bf 2)} The Restricted Isometry Property was proved in \cite{MPT1}
for matrices with independent rows (rather than columns), under a
sub-gaussian hypothesis.  It is worth noting that the corresponding
result for matrices with independent isotropic sub-gaussian columns
is not true in general.  One can see it by considering the matrix
with columns $X_i =
\sqrt{2}\delta_i(\varepsilon_{1i},\ldots,\varepsilon_{ni})$, where
$\delta_i$ are independent random variables, $\PP(\delta_i = 1) =
\PP(\delta_i = 0) = 1/2$ and $\varepsilon_{ji}$ are independent
Bernoulli variables, independent of $\delta_i$'s.  The vectors $X_i$
are then isotropic and sub-gaussian, but $\PP(X_i = 0) = 1/2$.  As a
consequence, the concentration hypothesis and thus the RIP property
are not satisfied.

\section{Main technical result}
\label{section:main technical result}

In this Section, $X_1,\ldots,X_N \in \R^n$ are independent $\psi_r$
random vectors for some (fixed) $0<r\leq 2$.
Let $1 \le m \le N$. We shall consider three quantities $A_m$, $B_m$
and $C_m$ depending on $X_1, \ldots, X_N$. Recall that $B_m$ has been
defined in (\ref{B_m_definition}) as
\begin{displaymath}
  \EE_m= \sup_{z\in U_m}\left| \left|\sum_{i\le N}z_iX_i
  \right|^2-\sum_{i\le  N}z_i^2|X_i|^2
 \right| ^{1/2}
\end{displaymath}
and define the other two quantities as follows:
$$
  \am = \sup_{z\in U_m}\left| \sum_{i \le N} z_i X_i \right| , \quad
\cm = \max_{i\leq
  N}  |X_i|.
$$
We clearly have
$$
   |\am^2 - \bm^2 | \leq \cm^2 .
$$
Given a real number $s$, we will denote $\max(s, 0)$ by $s_+$.


\medskip

The main purpose of this Section is
to prove Theorem \ref{UUP_thm}. In fact we will prove a stronger
technical result, Theorem \ref{UUPmain}, from which
Theorem \ref{UUP_thm} will follow.

\begin{Theorem}
\label{UUPmain}
Let $0<r\leq 2$. Let $n\geq 1$ and $1\le m\leq N$ be integers.
Let  $X_1,\ldots,X_N\in \R^n$ be independent $\psi_r$ vectors with
$\psi=\max_{i \le N}\|X_i\|_{\psi_r}$.
For every $1\leq m\leq N$,
$\theta\in (0,1/4)$, and  $K\geq 1$ one has
\begin{eqnarray}
\label{in-particular}
  \lefteqn{  \PP \left( B^2_m \geq \max\{B^2, \cm B, 24\, \th \,\cm ^2
      \}\r)} \nonumber\\
&\le &  (1+ 3\log m) \exp\left(- 2 K^r m^{(1+s)/2} \log \frac{2 N}{\th
    m}\r),
\end{eqnarray}
with $s=(1-r)_+$ and
$$
  B =  C_0^{1/r} \psi K \ m^{q-1/2} \left( \log \frac{2 N}{\th m}
    \r)^{1/r},
$$
where $C_0$ is an absolute constant  and $q=\max\{1, 1/r\}$.
\end{Theorem}

\medskip

\noindent {\bf Remark:}
In fact we shall prove a stronger statement:
with the notation of Theorem \ref{UUPmain},
for every $1\leq m\leq N$,
$\theta\in (0,1/4)$, and  $K\geq 1$,
and for  every $0\leq \ell \leq \log_2 m$,
 one has
 \begin{eqnarray}
   \label{general-ell}
\lefteqn{   \PP \left( B^2_m \geq \max\{\overline B^2, \cm \overline
    B, 24\, \th\, \cm ^2 \}\r) }\nonumber \\
   & \leq &  (1+ 2\ell)\exp\left(- 2 K^p \frac{m}{2^{\ell}} \log
  \frac{12 e N 2^{\ell}}{\th m}\r),
 \end{eqnarray}
where
$$
  \overline B = C^{1/r} \psi K \left( \left( \frac{m}{2^{\ell}} \r)^q
  \left( \log \frac{2 N 2^{\ell}}{\th m}  \r)^{1/r} + m^{q-1/2}
  \left( \log \frac{2 N}{\th m}  \r)^{1/r} \r),
$$
$C$ is an absolute constant and $q=\max\{1, 1/r\}$.

\bigskip

Before starting the proof of the theorem we show how  it implies
Theorem \ref{UUP_thm}, stated in Section \ref{section:RIP}.

\subsection{Proof of Theorem \ref{UUP_thm}}

Fix $K_1\geq 1$ and let $K\geq K_1$ be such that
$$
  K^2  m\log^{2/r} \frac{2N}{\theta m} = K_1^2 \theta^2 n.
$$
By Theorem~\ref{UUPmain} with $r\geq 1$, and
the condition on $m$,
\begin{eqnarray*}
  \lefteqn{
    \PP \left( B^2_m \geq \max\{B^2, \cm B, 24 \th \cm ^2\}\r)
}\\
   & \leq & (1+ 3\log m) \exp\left(- 2 K^r \sqrt m\log \frac{2
       N}{\th m}\r)\\
   & \leq &     \exp\left(- c K_1^r\sqrt m\log \frac{2 N}{\th m}\r)  ,
\end{eqnarray*}
where
$$
  B=  C_0 \psi K \, \sqrt{m}  \log^{1/r}\frac{2 N}{\th m} = C_0 \psi K_1
  \th \sqrt{n},
$$
and $c$ and  $C_0$ are absolute positive constants.
Thus, if $\cm \leq K_2 \sqrt{n}$  for some $K_2$, then
\begin{eqnarray*}
  \max\{B^2, \cm B, 24 \th \cm ^2\} & \leq & C_1 \th n
  \max\{\psi ^2 K_1^2, \psi K_1 K_2, K_2^2 \} \\
  & \leq &  C_1 \th n  \left( \psi  K_1 +  K_2 \r)^2 ,
\end{eqnarray*}
where $C_1$ is an absolute constant.
This concludes the proof.
\qed

\subsection{Proof of Theorem \ref{UUPmain}}

We will prove the theorem in a stronger form (\ref{general-ell}).
Then  (\ref{in-particular})  follows  by choosing $0\leq \ell\leq \log_2 m$
to be the largest integer satisfying
$$
  \frac{1}{2^{q \ell}} \left( \log \frac{2 N 2^{\ell}}{\th m} \r)^{1/r}
  \geq m^{-1/2} \left( \log \frac{2 N}{\th m}  \r)^{1/r} .
$$

The proof will use the same construction as in \cite{ALPT}, which
however requires some modifications. For completeness and the reader's
convenience we provide  details of the argument.

We require the following two lemmas proved in \cite{ALPT} with $r=1$.
Since the proofs for general $r$ repeat the same arguments, we leave
them for the reader.

\begin{lemma}
   \label{mainl}
Let $0<r\leq 2$ and $X_1,\ldots,X_N\in \R^n$ be independent $\psi_r$
vectors with $\psi=\max_{i \le N}\|X_i\|_{\psi_r}$. Let $m\leq N$,
$\eps, \alpha \in (0, 1]$. Let $q=\max\{1, 1/r\}$ and
$L\geq  m^q \left(2 \log \frac{12 e N}{m \eps}\r)^{1/r}$. Then
$$
 \PP \left(\sup _{F\subset \{1, ..., N\} \atop |F|\leq m } \ \sup
 _{E\subset F} \ \ \sup _{z\in {\cal{N}}(F, \eps, \alpha)} \
 \sum_{i\in E} \left| \la z_i X_i, \sum_{j\in F\setminus E} z_j X_j
 \ra \r| \geq \psi\, \alpha L A_m \r)
$$
$$
  \leq \exp\left(-\frac{1}{2} \ L^r \ m^{-(r-1)_+}\r).
$$
\end{lemma}

\begin{lemma}
\label{ltwo}
Let $0<r \leq 2$ and $X_1,\ldots,X_N\in \R^n$ be independent $\psi_r$
vectors with $\psi=\max_{i \le N}\|X_i\|_{\psi_r}$.
Let $1\leq k,  m\leq N$, $\eps, \alpha \in (0, 1]$, $\beta >0$, and $L>0$.
Let $B(m, \beta)$ denote the set of vectors $x\in \beta B_2^N$ with
$|\supp  x| \leq m$ and let $\cal{B}$ be a subset of $B(m, \beta)$ of
cardinality $M$. Then
\begin{eqnarray*}
\PP \left(\sup _{F\subset \{1, \ldots, N\} \atop |F|\leq k }\ \sup
  _{x\in {\cal{B}}} \ \ \sup _{z\in {\cal{N}}(F, \eps, \alpha)} \right.
 & &\!\!\!\!\!\! \left. \sum_{i\in F} \left| \la z_i X_i,
     \sum_{j\not\in F}    x_j X_j \ra \r|
    \geq \psi \alpha \beta L A_m \r)  \\
  &\leq &  M \left( \frac{6 e N}{k \eps} \r)^k
  \exp\left(-\frac{1}{2} \ L^r \  k^{-(r-1)_+}\r).
\end{eqnarray*}
\end{lemma}

The following formula is well known and the proof is in its statement.

\begin{lemma}\label{dec}
Let $x_1, \ldots, x_N \in \Rn$, then
$$
 \sum_{i\ne j} \la x_i, x_j \ra  =4 \cdot 2^{-N}
  \sum_{E\subset \{1, ..., N\} }\sum_{i\in E}
  \sum_{j\in E^c} \la x_i, x_j \ra.
  $$
\end{lemma}

\medskip

We are now ready to start the proof of Theorem~\ref{UUPmain}.

\noindent {\bf Proof of Theorem~\ref{UUPmain}.}
As in \cite{ALPT}, the construction   splits into two cases.

If $\ell=0$
we set
$$
   {\cal{M}(\theta)} = \bigcup _{E\subset \{1,\ldots N\} \atop |E|=m}
   {\cal{N}}(E, \theta/4, 1).
$$
Otherwise,
 define positive integers $a_0, a_1,\ldots, a_\ell$ by
$a_k:= [m\, 2^{-k+1}]-[m\,
2^{-k}]$ for $1 \le k \le \ell$ and $a_0:= [m\, 2^{-\ell}\,]$.
Observe that
$a_k \le m\, 2^{-k+1}$ for $ 1 \le k \le \ell $, $a_0 \le m\,2^{-\ell}$
and $\sum_{k=0}^\ell a_k = m.$
Recall that for $E\subset\{1,\dots, N\}$ we identify $\R^E$ with the
subspace of vectors in $\R^N$
with coordinates supported by $E$.

We consider $(\ell +1)$-tuples
$\left( (E_0, x_0), \dots, (E_\ell,x_\ell)\right)$
where
$(E_k)_{0\le k\le \ell}$  are mutually disjoint subsets of
$\{1,\ldots N\}$, $\  |E_k|\le a_k$, $\  x_k\in\R^{E_k}$ for all
$0\le k\le \ell$.
A $(\ell +1)$-tuple
$\left((E_0, x_0), \dots, (E_\ell,x_\ell)\right)$ is  said to be
admissible if
$$
 x_k\in {\cal{N}}\left(E_k, \theta 2^{-k}, \sqrt{{2^k\over m}}\right)
   \mbox{for \,} 1\le k\le \ell,
 x_0\in {\cal{N}}\left(E, \theta/4,1\right),
\left|\sum_{k=0}^{\ell}x_k\right|\le 2.
$$

The set of all vectors $x=\sum_{k=0}^{\ell}x_k$ associated to
admissible $(\ell +1)$-tuples $\left((E_0, x_0), \dots,
  (E_\ell,x_\ell)\right)$ will be denoted by ${{\cal{M}}(\theta)}$.

We shall consider the details of the case
$\ell >0$,
the other case
can be treated similarly.

Fix $\left((F_0, x_0), \dots, (F_\ell,x_\ell)\right)$ to be admissible
and let $x= \sum _{k=0}^{\ell} x_k \in {\cal{M}}(\theta)$.
Denote the coordinates of $x$ by $x(i)$,
$i\leq N$, then
\begin{equation*}
  |A x|^2  =\la \sum _{i\leq N} x(i) X_i, \sum _{i\leq N} x(i) X_i \ra
  =   \sum _{i\leq N} x(i)^2 |X_i|^2 + \sum _{i\ne j} \la x(i) X_i,
  x(j) X_j  \ra.
\end{equation*}
So
\begin{equation}
\left| |A x|^2 - \sum _{i\leq N} x(i)^2 |X_i|^2\right|=\left| D_x\right|
\label{decompos}
\end{equation}
where
$$
     D_x=\sum _{i\ne j} \la x(i) X_i, x(j) X_j \ra.
$$
Now we split $D_x$ according to the structure of $x$. Namely
we let
$$
   D'_x := \sum _{k=0}^\ell \sum _{i, j \in F_k \atop i\ne j} \la x(i)
   X_i, x(j) X_j \ra,
$$
and
$$
   D''_x : =  \sum _{k=0}^\ell \sum _{i \in F_k \atop j\not\in F_k} \la
   x(i) X_i, x(j) X_j \ra ,
$$
so that we have
$$
 \left| |A x|^2-   \sum _{i\leq N} x(i)^2 |X_i|^2 \right|=
     \left| D'_x+ D''_x\right| \le \left| D'_x\right|+\left|
     D''_x\right|.
$$

We first estimate $D'_x$.
By Lemma~\ref{dec} we have
\begin{eqnarray*}
   D'_x &=& \sum _{k=0}^\ell \sum _{i, j \in F_k \atop i\ne j} \la x(i)
   X_i, x(j) X_j \ra   \\
& = &
4\sum _{k=0}^\ell
2^{-|F_k|}   \sum_{E\subset F_k }\sum_{i\in E}
  \sum_{j\in F_k\backslash E} \la x(i) X_i, x(j) X_j   \ra.
\end{eqnarray*}
Thus
$$
\left| D'_x\right|
\le
4 \sum _{k=0}^\ell
2^{-|F_k|}   \sum_{E\subset F_k }\left|
\sum_{i\in E}
  \sum_{j\in F_k\backslash E} \la x(i) X_i, x(j) X_j   \ra\right|
$$
$$
  \le 4\sum _{k=0}^\ell
\sup_{E\subset F_k }\left| \sum_{i\in E}
  \sum_{j\in F_k \backslash E} \la x(i) X_i, x(j) X_j   \ra\right|
$$
and  using the fact that
  $|F_k|\le a_k $
for $ 0 \le k \le \ell$,
we arrive at
$$
 \left| D'_x\right|\le
   4\sum _{k=0}^\ell
\sup_{F\subset \{1, ..., N\} \atop |F|\leq a_k}\
\sup_{E\subset F }\sum_{i\in E}
\left|  \la x(i) X_i, \sum_{j\in F\backslash E} x(j) X_j   \ra\right|.
$$

We now set  $q=\max\{1, 1/r\}$ and apply Lemma~\ref{mainl} to each summand
in the sum above with the parameters
$$
  a_0, \eps = \theta/4, \alpha =1, \mbox{\ and }
  L= K \left( \frac{m}{2^\ell} \r) ^{q} \left( 4
  \log\frac{ 48 e N 2^\ell}{\theta m} \r) ^{1/r}
$$
for $k=0$ and
$$
  a_k, \eps =\theta  2^{-k}, \alpha = \sqrt{\frac{2^k}{ m}}, \mbox{\ and }
  L= K \left( \frac{m}{2^k} \r) ^{q} \left( 4 \log\frac{ 12 e N
  4^k}{\theta m} \r) ^{1/r}
$$
for $ 1 \le k \le \ell$.
By the union bound we obtain
that the probability of the event
\begin{eqnarray*}
     \sup _{x\in {\cal{M}}(\theta)} |D'_x|  &\geq &
  \psi A_m K \left( \left( \frac{m}{2^\ell} \r) ^{q} \left( 4
  \log\frac{ 48 e N 2^\ell}{\theta m} \r) ^{1/r} \right.\\
    & + &  \left.\sum _{k=1}^\ell \left( \frac{m}{2^k} \r) ^{q-1/2}
  \left( 4 \log\frac{ 12 e N 4^k}{\theta m} \r) ^{1/r} \r)
\end{eqnarray*}
is not larger than
$$
    \exp\left( - K^r \ \frac{2m}{2^\ell} \log\frac{48 e N
        2^{\ell}}{\theta m} \right)
   +  \sum _{k=1}^\ell
   \exp\left( - K^r \ \frac{2m}{2^k}  \log\frac{12 e N 4^k}{\theta m}
   \right) .
$$
Therefore the probability of the event
\begin{eqnarray*}
   \sup _{x\in {\cal{M}}(\theta)} |D'_x|  & \geq &
  \psi A_m K \left( \left( \frac{m}{2^\ell} \r) ^{q} \left( 4
  \log\frac{ 48 e N 2^\ell}{\theta m} \r) ^{1/r} \right.\\
   & + & \left.  C_1^{1/r} m^{q-1/2} \left( \log\frac{2 N}{\theta m}
       \r) ^{1/r} \r)
\end{eqnarray*}
is not larger than
$$
   (1+\ell) \exp\left( - K^r \ \frac{2m}{2^\ell}
  \log\frac{12 e N 2^{\ell}}{\theta m} \r) ,
$$
where $C_1$ is an absolute constant.

We now pass to the estimate for $D''_x$ which essentially follows the
same lines.

For every $1\leq k \leq \ell$ we consider ${\cal{M}}_k(\theta) =
{\cal{M}}_k'(\theta) \cap 2B_2^N$, where ${\cal{M}}_k'(\theta)$
consists of all vectors of the form $v= v_0 + \sum _{s=k+1}^{\ell}
v_s$, where $v_i$'s ($i=0, k=1, \ldots, \ell$) have pairwise disjoint
supports and
$$
v_0 \in \bigcup _{E\subset \{1,\ldots N\} \atop |E|\leq a_0}
{\cal{N}}(E,\theta/4, 1) , \, \, v_s \in \bigcup _{E\subset \{1,\ldots
  N\} \atop |E|\leq a_s} {\cal{N}}\left(E, \theta\, 2^{-s},
  \sqrt{\frac{2^s}{m}}\r) \, \mbox{ for } \, s \ge k+1 .
$$
Then ${\cal{M}}_k (\theta) \subset 2B_2^N$ and
(similarly as in   \cite{ALPT})
we can estimate the cardinality
\begin{eqnarray*}
  |{\cal{M}}_k (\theta)| &\leq &  \left({12\over \theta}\right)^{a_0}
  \prod _{s=k+1}^\ell \left({3\cdot     2^s\over\theta}\r)^{a_s}
  {N \choose a_s} \leq \left({12\over \theta}\right)^{a_0}
   \prod _{s=k+1}^\ell \left(\frac{3\cdot 2^s  e N}{\theta a_s}\r)^{ a_s}   \\
   & \leq &  \exp\left( \sum _{s=k+1}^{\ell+1} \frac{2m}{2^s} \log
     \frac{3e 4^s        N}{2\theta m} \r)   \leq
\exp\left( \frac{4 m}{2^k} \log \frac{6 e 4^k N}{\theta m}
     \r) .
\end{eqnarray*}

Recalling that $x= \sum _{k=0}^{\ell} x_k \in {\cal{M}}(\theta)$ for
some admissible $(\ell +1)$-tuple  $\left((F_0, x_0), \dots,
(F_\ell,x_\ell)\right)$  and
setting $G_k = \{0, k+1, k+2, \ldots , \ell\}$, we observe that
\begin{eqnarray*}
  |D''_x |  & = &  \Big|2 \sum _{k=1}^\ell  \sum _{i \in F_k } \la
  x(i) X_i, \sum
  _{r\in G_k}  \sum _{j\in F_r}   x(j) X_j \ra \Big|\\
 & \leq &  2 \sum _{k=1}^\ell  \sup _{F\subset \{1, ..., N\} \atop
    |F|\leq 2m/2^k}\
  \sup _{u\in {\cal{N}}(F, 2^{-k}, \sqrt{2^k/m})} \sup _{ v\in
    {\cal{M}} _k }
  \sum_{i\in F} \left| \la u(i) X_i, \sum_{j \not\in F} v(j) X_j \ra \r|  .
\end{eqnarray*}

Now we apply Lemma~\ref{ltwo} to each summand
$k=1, \ldots, \ell$,
with parameters
$$
 \eps =\theta 2^{-k}, \alpha =\sqrt{2^k\over m},  \beta = 2,
  {\cal{B}}  ={\cal{M}} _k(\theta) \mbox{\ and \ }
  L = K \left( \frac{m}{2^k} \r) ^{q} \left( 12 \log\frac{ 12 e N
  4^k}{\theta m} \r) ^{1/r}.
$$
Using the union bound we obtain
\begin{eqnarray*}
\lefteqn{  \PP \left( |D''_x| \geq 2  \psi A_m K \sum _{k=1}^\ell
   \left( \frac{m}{2^k} \r)^{q-1/2}  \left(12
   \log \frac{12 e N 4^k}{\theta m} \r)^{1/r} \r) } \\
&\leq &
  \sum _{k=1}^\ell \exp\left(\frac{4 m}{2^k} \log \frac{12 e 4^k
    N}{\theta m} +
  \frac{2 m}{2^k} \log \frac{3 e 4^k N}{\theta m} -
  K^r\ \frac{12 m}{2^k} \log \frac{12 e N 4^k}{\theta m}\r) \\
 & \leq &  \sum _{k=1}^\ell  \exp\left( - K^r\ \frac{6 m}{2^k} \log
  \frac{12 e N 4^k}{\theta m}\r) \leq \ell \exp\left( - K^r \
  \frac{6 m}{2^\ell} \log \frac{12 e N 4^\ell}{\theta m}\r).
\end{eqnarray*}
Thus
\begin{eqnarray*}
  \lefteqn{
  \PP \left( \sup _{x\in {\cal{M}}(\theta)} |D''_x| \geq C_2^{1/r} \psi
  A_m  K  m^{q-1/2}    \left(\log \frac{2 N}{\theta m} \r)^{1/r}  \r) }\\
  & \leq &  \ell \exp\left( - K^r \
  \frac{6 m}{2^\ell} \log \frac{12 e N 4^\ell}{\theta m}\r) ,
\end{eqnarray*}
where $C_2$ is the an absolute constant.

\medskip
Since $D_x = D'_x + D''_x$, then
\begin{equation} \label{d-x-estimate-0}
  \PP \left( \sup _{x\in {\cal{M}}(\theta)}  |D_x| \geq A_m \gamma \r)
  \leq \left( 1 + 2 \ell \r) \exp\left( - K^r \
  \frac{2 m}{2^\ell} \log \frac{12 e N 2^\ell}{\theta m}\r) ,
\end{equation}
where
$$
  \gamma =  C_3^{1/r} \psi K \left( \left( \frac{m}{2^\ell} \r) ^{q} \left(
  \log\frac{ 2 N 2^\ell}{\theta m} \r) ^{1/r}
   + m^{q-1/2} \left( \log\frac{2 N}{\theta m} \r) ^{1/r} \r)
$$
for some absolute constant $C_3>0$.

\bigskip

Passing now to the approximation argument,
pick an arbitrary  $z \in S^{N-1}$  with $|\supp z|\leq m$.
 Define the following subsets of $\{1, \ldots, N\}$ depending on $z$.
Denote the coordinates of $z$ by $z(i)$ ($i=1, \ldots, N$).
Let  $n_1, \ldots, n_N$ be such that
$|z(n_1)|\geq |z(n_2)| \geq \ldots \geq |z(n_N)|$, so that
$z(n_i) = 0$ for $i>m$ (since $|\supp z|\le m$).
If $\ell=0$,
we denote the support of $z$ by $\widetilde E_0$ and consider only
this $\widetilde E_0$.
Otherwise  we set
$$
  \widetilde E_0 = \{n_i\}_{1\leq i \leq m/2^{\ell}}
$$
and
$$
  \widetilde E_1 = \{n_i\}_{m/2 < i\leq m}, \ \widetilde  E_2 =
  \{n_i\}_{m/4 < i\leq m/2}, \
  \ldots, \ \widetilde  E_\ell = \{n_i\}_{m/2^\ell < i \leq  m/2^{\ell-1}}.
$$
Recall that  integers $a_0, a_1,\ldots, a_\ell$ have been defined at
the beginning of this proof. Then, clearly,
$$
  a_0  = |\widetilde  E_0| \leq m/2^\ell, \quad
 a_k = |\widetilde  E_k|\leq m/2^k + 1 \leq m/2^{k-1} \, \,  \mbox{
   for every }   1\leq k\leq \ell,
$$
and $\sum _{i=0}^\ell a_i=m$.
 Also observe
that, since $z \in S^{N-1}$, then  for every $k\geq 1$,
$$
    \| P_{\widetilde  E_k} z \| _{\infty} \leq |z(n_s)| \le
    \sqrt{\frac{2^k}{m}},
$$
where $s = [m /2^{k}]$.

For every $k\geq 1$ the vector $P_{\widetilde E_k} z$ can be
approximated by a vector from ${\cal{N}}\left(\widetilde E_k,\theta
  2^{-k}, \sqrt{\frac{2^k}{m}}\r)$ and the vector $P_{\widetilde E_0}
  z$ can be approximated by a vector from ${\cal{N}}(\widetilde E_0,
  \theta/4, 1)$.  Thus there exists $x\in {\cal{M}}(\theta)$, with a
  suitable representation $x = \sum_{k=0}^\ell x_k$, such that
$$
   |z-x|^2 \leq \sum _{k=0}^\ell
|P_{\widetilde  E_k} z - x_k |^2 \leq
  \theta^2( 2^{-4} + \sum _{k=1}^\ell 2^{-2k}) < \theta^2\,(0.4).
$$
Moreover, $x$ is chosen to have the same support as $z$, and thus
$w = z-x$  has the support $|\supp w | \le m$.

It follows from the definitions of $D_z$ and $A$  that
$$
D_z = D_x + \langle Aw, Ax \rangle  + \langle Az, Aw \rangle
- \sum_{i\le N} w(i) \left(x(i) + z(i)\right) |X_i|^2,
$$
(here $w(i)$, $x(i)$
and $z(i)$
 denote the coordinates of $w$, $x$ and $z$, respectively).
Thus
$$
|D_z|\le |D_x|+ |Aw| (|Ax|+ |Az|) + |w|\,|x+z|\,\max_{i\le N}
|X_i|^2.
$$
It follows that
$$
   B_m^2 = \sup _{z\in S^{N-1} \atop |\supp z| \leq m} |D_z| \leq
   \sup _{x \in {\cal{M}}(\theta)} |D_x|
   +  2 \theta \left(A_m^2 + \cm ^2 \r) \leq
   \sup _{x \in {\cal{M}}(\theta)} |D_x|
   +  2 \theta \left(B_m^2 + 2 \cm ^2 \r) .
$$
Thus, by (\ref{d-x-estimate-0}) and using again $\am \leq
\sqrt{\bm ^2 + \cm ^2} \leq \bm +\cm$ we obtain
$$
  \PP \left( (1-2 \th) B_m^2 \geq 4 \th \cm ^2 + \cm \gamma + \bm \gamma
  \r) \leq \left( 1 + 2 \ell \r) \exp\left( - K^r \ \frac{2 m}{2^\ell} \log
  \frac{12 e N 2^\ell}{\theta m}\r) .
$$
Since $\th \leq 1/4$,  this  implies
$$
  \PP \left( B_m^2 \geq \max\{ 24 \th \cm ^2, 6 \cm \gamma, 6 \gamma^2 \}
  \r) \leq \left( 1 + 2 \ell \r) \exp\left( - K^r \ \frac{2 m}{2^\ell} \log
  \frac{12 e N 2^\ell}{\theta m}\r) ,
$$
which completes the proof.
\qed

\subsection{Optimality of estimates}

We conclude this section by an example showing optimality, in a
certain sense,  of
estimates in Theorem \ref{UUP_thm}.
We will limit ourselves  to the $\psi_1$ case, that is to  $r=1$.
%
To this end we consider a special case when $X_i = (X_{ij})_{j=1}^n$
where $X_{ij}$ are i.i.d. symmetric exponential variables with
variance one. We begin by showing an optimal estimate for $A_m$.

 First, from  \cite{ALPT} (Theorem 3.5) we have that for $N \le
\exp(c\sqrt{n})$ and any $K\geq 1$,
\begin{equation}
\label{A_m estimat}
  \PP \left( A_m \geq C K \left(\sqrt{n} + \sqrt{m}
  \log \frac{2N}{m}  \r)  \r)  \leq \exp\left( - c K \sqrt{n} \r)
\end{equation}
where $C,c>0$ are numerical constants. In the other direction, we
have the following
\begin{prop}
  \label{exponential}
  For any $1 \le m \le N$ and $t \ge 1$,
\begin{displaymath}
\PP\Big(A_m \ge c\Big(\sqrt{n} + \sqrt{m}\log\Big(\frac{2N}{m}\Big) +
    t \Big) \Big) \ge c\wedge e^{-t},
\end{displaymath}
where $c > 0$ is an absolute constant.
\end{prop}

Before we prove this proposition let us explain its relevance to
Theorem \ref{UUP_thm}.  Firstly, observe that the proposition implies
that with probability bounded away from zero, $A_m \ge c(\sqrt{n} +
\sqrt{m}\log(2N/m))$.  This shows that -- except for allowing a change
of absolute constants -- one cannot obtain a better bound on $A_m$
than (\ref{A_m estimat}), valid with overwhelming probability (i.e.,
with probability converging to one as $n \to \infty$).
Secondly, assume that $N \le \exp(c\sqrt{n})$.  By taking $t = c K
\sqrt{n}$, we obtain that for large $n$, $\PP(A_m \ge cK\sqrt{n}) \ge
\exp(-cK\sqrt{n})$.  We compare this with estimates for
probabilities in (\ref{A_m estimat}). Namely, using
Lemma \ref{max} (noting that the density of $X_i$'s is log-concave),
we  can see that for $m \log^2(2N/m) \le n$, the theorem
implies that $\PP(A_m \ge CK\sqrt{n}) \le
\exp(-\tilde{c}K\sqrt{n})$. So in this range of $m$ the upper and
lower bounds on probability coincide up to numerical constants in the
exponent.

Regarding Theorem \ref{UUP_thm}, again assume that $N \le
\exp(c\sqrt{n})$. Using again Lemma \ref{max}, we get with
overwhelming probability that for all $i$, $|X_i| \le C'\sqrt n$.
Now assume that for some $m$ we have with overwhelming probability
$B^2_m \le Cn$.  Then by the obvious bound $A_m^2 \le B^2_m +
\sup_{z \in
  U_m}\sum_{i\le N}|z_i|^2 |X_i|^2$, with probability close to one we
also have $A_m \le C''\sqrt{n}$.  On the other hand, as noted above,
$\PP(A_m \ge c(\sqrt{n} + \sqrt{m}\log(2N/m)))$ is bounded away from
zero. Thus,
$ c(\sqrt{n} + \sqrt{m}\log(2N/m)) \le C''\sqrt n$, which in turn
implies that for $n$ large enough we have $m\log^2(2N/m) \le Cn$.
This shows that the factor $\log^2(2N/\theta m)$ in Theorem
\ref{UUP_thm} is of the right order.

\vskip0.5cm

\noindent{\bf Proof of Proposition \ref{exponential}{\ \ }}
Since
\begin{displaymath}
  A_m = \sup_{\alpha \in S^{N-1}\atop |\supp \alpha| \le m}\sup_{\beta
    \in S^{n-1}}\sum_{ij}\alpha_i \beta_j X_{ij},
\end{displaymath}
by general tail estimates for linear combinations of exponential
variables with vector valued coefficients (see e.g. Corollary 1 in
\cite{L}), we get
\begin{displaymath}
\PP\Big(A_m \ge c\big(\E A_m + \sqrt{t}\sigma + t b \big) \Big) \ge
c\wedge e^{-t},
\end{displaymath}
where
\begin{displaymath}
\sigma^2 = \sup_{\alpha \in S^{N-1}\atop |\supp \alpha| \le
  m}\sup_{\beta \in S^{n-1}}\sum_{ij}\alpha_i^2\beta_j^2 = 1
\end{displaymath}
 and
\begin{displaymath}
b = \sup_{\alpha \in S^{N-1}\atop |\supp \alpha| \le m}\sup_{\beta \in
  S^{n-1}}\max_{ij}|\alpha_i\beta_j|=1.
\end{displaymath}

Therefore, it is enough to show that $\E A_m \ge c(\sqrt{n} +
\sqrt{m}\log(2N/m))$. Obviously, $\E A_m \ge c\sqrt{n}$, since a
single column of the matrix $A$ has expected Euclidean norm of the
order $\sqrt{n}$. As for the other term, it is enough to consider the
first row of our matrix.  We have
\begin{displaymath}
\sqrt{m} A_m \ge \sup_{\alpha \in \{0,-1,+1\}^N\atop |\supp \alpha| =
  m} \sum_{i=1}^N \alpha_i Y_i,
\end{displaymath}
where to simplify the notation we set $Y_i = X_{i1}$. On the right
hand side we actually have $\sum_{i=1}^m |Y_i^\ast|$, where $Y_i^\ast$
is such a rearrangement of $Y_i$ that $|Y_1^\ast|\ge |Y_2^\ast| \ge
\ldots \ge |Y_n^\ast|$, which can be used to derive lower bounds on
the expectation. We will however not rely on this representation,
instead we will use a Sudakov type minoration principle for
exponential variables proved in
\cite{T1}, Theorem 5.2.9, which we state here in a simplified version,
adapted to our purposes.

\begin{lemma}
Let $Y_1,\ldots,Y_N$ be independent symmetric
  exponential variables with variance one. Consider $T\subseteq
  \ell_2^N$ of cardinality $k$ and $u \ge 1$. If for any $s,t \in T$,
  $t\neq s$,
\begin{displaymath}
\sqrt{u}|t-s| + u\|t-s\|_\infty > u,
\end{displaymath}
then $\E \max_{t\in T} \sum_{i=1}^N t_iY_i \ge c\min(u,\log k)$, where
$c > 0$ is an absolute constant.
\end{lemma}

In our case, $T = \{\alpha \in \{0,-1,1\}^N\colon |\supp \alpha|\le
m\}$, so $k \ge \binom{N}{m}$. Also, since $\|t-s\|_\infty \ge 1$
for $t,s \in T$, $t\neq s$, the condition of the lemma is trivially
satisfied for any $u \ge 1$, in particular for $u = \log k$. Thus,
for $m \le N/2$, we obtain $\sqrt{m}\,\E A_m \ge \log k \ge cm\log
(2N/m)$. On the other hand we have $\E A_m \ge c\sqrt{m}$, so for $m
\ge N/2$ it is enough to adjust the constants.  \qed

\vspace{1cm}

\address

\end{document}